# Productivity equation and the *m* distributions of information processing in workflows

Charles Roberto Telles

**Abstract**

This research investigates an equation of productivity for workflows regarding its robustness towards the definition of workflows as probabilistic distributions. The equation was formulated across its derivations through a theoretical framework about information theory, probabilities and complex adaptive systems. By defining the productivity equation for organism-object interactions, workflows mathematical derivations can be predicted and monitored without strict empirical methods and allows workflow flexibility for organism-object environments.

## I. INTRODUCTION

Mostly structural organization emphasize methods considering only physical aspects of routines (discrete variables) such as products properties or services quality, both identifiable with several methodologies nowadays [1]. But for a hybrid organization, subjective work (agents) can't be definable by such methodologies and that aspect is one of the main points in this research. Nowadays general and famous international industries methods of workflows consider mostly discrete variables only as the basis of production. But several other types of industries or labor activities demand continuous variables to be defined within a method. Following this path, this article considers productivity in modern contexts of human interactions with non-physical environments [1] or in other words, hybrid structural organizations where non-physical information is the main component of productivity.

A large number of variables constitute workflows and are complex in nature [2] as a complex adaptive system [3]. Statistical analyzes on the portion that occupies each category of information processing of a workflow on a scale of 0 to 100% is not recommended for analyzing how an organism process information due to the nonlinear nature of the phenomena [4]. The nonlinearity of the complex adaptive systems prevents the possibility of comparing a group of individuals with different specific and indeterminate cognitive patterns naturally [5], and their forms of work, to generate patterns of execution in the workflows of linear or nonlinear dimensions [6-9].

## II. PRODUCTIVITY EQUATION

**Definition 1.** The ergodicity of a system as workflows [10-12], is not a constantly expression considering real-life situations. Therefore, for the purpose of empirical investigation of actual facts, it is recommended starting from an analysis in which non-ergodicity is the a priori event, more present in the real world, where the distribution *m* assumes various forms (derivations) [12-16].

If a stream of information arises from an event $i$, where individual $X$ processes as an input state a given discrete information $I$ generating time $T$ and reaching precision $P$ influenced by individual experience $I^i$, then there is a probabilistic event in which the event $i$, the probability of precision ($P$) can be defined as $Prob(P) = (I \text{ and } T)$ where $Prob\ X(I \text{ and } T) = [P = i_1, \dots, x_n = i_n] = prob_{i,\dots,i_n}$. Thou, as $T$ is defined by $I$ given (data source) and is processed by individual $X$ defining how information will be processed and how time will be generated as $Prob\left(\frac{T}{I}\right)$, the input values of $I$ and $I^i$ generates $m$ distributions for $I$ influenced by $I^i$ and for $T$ also influenced by $I^i$ as an output. $I^i$ can be defined as the *ad hoc* cognitive performances of the individual in which $I^i$ represents a posterior processing stage of $Prob\left(\frac{T}{I}\right)$, being this first probability the first interactions between organism and object. In this sense, $I^i$ is the same event $Prob\left(\frac{T}{I}\right)$, but exponentially growth by individual experience $i$, being $i$ in its turn defined by several others probabilities generated by the interactions and iterations of the organism and environment.

This equation of event $i$ as any event that considers information processing by an organism giving $P$, can be defined by,

$$i = Prob(I^i).Prob\left(\frac{T}{I}\right) = m, \tag{1}$$

In which $P$ assumes $m$ variations according to the probabilistic distribution of $I$, $I^i$ and $T$, which in turn defines the information entropy of $P$, hence, of the workflow (see $m$ distributions part III). Also the equation (1) can be written in the same way as equation (2) due to conceptual description mentioned before.

$$i = Prob\left(\frac{T^i}{I}\right).Prob\left(\frac{T}{I}\right) = m \text{ (7) or } i = Prob\left(\frac{T^{\frac{T}{I} \to \infty}}{I}\right).Prob\left(\frac{T}{I}\right) = m, \tag{2}$$

It is necessarily understood that the time occurs only if the information is processed, so the probability of time is dependent on the probability distributions of the information and the processing. However, the probabilistic distributions of $I$ and $T$ assume behavior in a sample space that does not have fixed intervals, since they come from complex adaptive systems [2-5] and with a degree of freedom for any resultant that varies from individual to individual [5]. In this way, it is possible to assume that every learning process as well cognitive processing derives not from a predefined sample understanding it as like information is fully objective for the organism in its potential of apprehension. For this reason, it can be identified as a probabilistic event defined by the data source ($I$), individual experience ($I^i$) as input values and time of processing ($T$) and resulting information ($I$) as output values. Thus, it is not previously defined that there are probabilities of $I$ and $T$, but that the probabilities are due to the dependence between these variables and their empirical expression, which identifies the phenomenon itself as the product of nonlinear dynamics by cumulatively modifying their probabilistic distributions as the more complex adaptive systems interact with the environment.

If theoretically an information processed in $t$ time is set to defined variables and processing in $I$, the precision is always given as 100%. But, if $I^i$ influence the processing of $I$, then $I^i$ can influence the event to express monotone properties where $I^i$ takes chance for lowering or rising $Prob\left(\frac{T}{I}\right)$ chances. Lowering means lower precision or time processing and rising means achieving new form of processing the formal way of dealing with $Prob\left(\frac{T}{I}\right)$. It means in other words, improving the $Prob\left(\frac{T}{I}\right)$ default pattern of working. Consider for this a single example of empirical experiments of eye processing forms [17], where statistical patterns are processed in the task of recognizing an individual's own error of visual memory and learning. As different tasks are presented by the environment, the default pattern of information processing adapts itself to reach the ideal precision of the given task. It means clearly that organism influences information and time processing performance. The same way, information as it is given influences the organism to adapt itself in order to achieve specific goals. For this experiment, chances of precision in the given events present lowering or rising effects. Another [18] research, even if subjectivity of an individual is removed, empirically, the response of information processing can be changed from the neural perspective between individuals and also with specific processing capacities.

## III. *M* DISTRIBUTIONS OF INFORMATION PROCESSING IN WORKFLOWS

The mathematical modelling that describe information flows in the workflows can be defined as a probability event when from a given event $i$ by the equation (1) of the interaction between given discrete information variables ($I$), individual experience $I^i$ and defining the time $T$ (discrete or continuous depending on the case) as a function of *I* for the execution of individual or intersectorial work/between agents reaching precision (*P*), assumes for *i* expressions like:

*A. Discrete binomial probability distribution*

**Definition 2.** For a single event and *i.d.d*, in which *P* assumes *m* variations according to the probabilistic distribution of *I*, $I^i$ and *T*, which in turn defines the information entropy of *P*, where *i* assumes expressions like:

$$if\ i = 1 \therefore fi = Prob(I^1).Prob\left(\frac{T}{I}\right) = Prob\ (T), \qquad (3)$$

Where *T* is the unique variable on system defined by successful *p*. Or, *if* $i = 0$,

$$i = Prob(I^0).Prob\left(\frac{T}{I}\right) = Prob\left(\frac{T}{I}\right), \qquad (4)$$

Where decision making takes place in general defining *T* and *I* as,

$$f(i;p) = \begin{cases} p & if\ i = 1 \\ q = 1 - p & if\ i = 0 \end{cases} \qquad (5)$$

*B. Discrete probability distribution*

**Definition 3.** Otherwise, if *m* be discrete, not binomial and information is processed constantly without imprecisions in the same properties as cognition processing, reaching 100% precision for any event,

$$\sum_i Prob(P = i) = i_1, \dots, x_n = i_n = p_{i,\dots,i_n} \qquad (6)$$

Where *i* is processed constantly without any oscillation, then,

$$\sum_{i_1=0}^{1} \cdots \sum_{i_n=0}^{1} p_{i_1,\dots,i_n} = 1 \qquad (7)$$

$$i = Prob(I^i = 1).Prob\left(\frac{T}{I}\right) = \left(\frac{T}{I}\right), \qquad (8)$$

Where $\left(\frac{T}{I}\right)$ assumes main role of processing as a defined precision given by organism. Whereas the precision is not reached at 100%, indefinite probabilities occurs due to subjective nature of organism processing, but as for an effort of reaching 100% precision, cumulative trials leads the event to the maximum information entropy where trials of learning process assumes also monotone properties as following items briefly describes.

*C. Discrete probability distribution and monotonically decreasing*

Defined by equation $F_p(I) = Prob(P > I)$ from a discrete cumulative distribution function (CDF) in which the discrete and probabilistic function is discrete by the presence of pre-defined data, but with uncertain processing and / or temporality that assume in an adaptive system the well-defined *ad hoc* mode of work and predictable. This type of distribution associated with adaptive systems does not represent an entropy of critical information, in which it is not possible to reach accuracies close to 100%. The monotonically decreasing function is presented by the high probability of accurate execution and low presence of information that generates randomness in the system. It is dependent on an *ad hoc* method to achieve an accuracy of 100 or close. Mathematically, it is not possible to obtain a CDF-like probability as described above as a function of the organic component of the system. Thus, the cumulative function of information as discrete values can of course be processed to the inverse of the manifestation itself in its physical nature or axiomatic origin of probabilities (a characteristic event of a non-adaptive system where $Prob(P < I)$). Precisely the propositions of biological order can establish a function between precision and information of the type $F(P) > F(I)$. Thus the probability of precision is strictly greater than that of information, the set of $P$ contained in $I$ being at the same time contained in another set of unknown dimension (of individual experience), such as, $\mathbb{P}((\infty, P]) \leq \mathbb{P}((-\infty, I]) \leq \mathbb{P}((-\infty, I^i])$, being $I^i$ the experience accumulated by the individual, or in other words, the accumulated information of $n$ events $i$, which confer to the biological potential the possibility that $\mathbb{P}((\infty, I^i] > \mathbb{P}((-\infty, P]))$. A cumulative function in an adaptive system assumes the biological form of the individual and breaks the axiom of probabilities, differentiating the axiom that applies to the physical world from the complex adaptive world.

**Definition 4.** In the other hand if $I^i > I$, for any $I$ given to prevent oscillations on system resulting in $F_p(I) = Prob(P > I)$, thus the probability of precision is strictly greater than that of information, the set of $P$ contained in $I$ being at the same time contained in another set of unknown dimension (of individual experience $(I^i > I)$), such as,

$$i = Prob(I^i > I).Prob\left(\frac{T}{I}\right) = Prob\left(\frac{T}{I}\right)^{I^i > I} \tag{9}$$

Which confer to the biological potential the possibility,

$$\mathbb{P}((\infty, I^i] > \mathbb{P}((-\infty, P])). \tag{10}$$

To add biological properties in the complex adaptive system, learning process can be observed as a heuristic input and output of information of obscure probabilities. But for any given heuristic cognitive processing, weights can be strictly associated with organism search for environment patterns and previous memory experiences.

Weights in this view are distributed and classified as cognitive modeling of cumulative experiences $I^i$ as cognitive processing are indexed heuristically as $1, \ldots, n$ where $I_1^i, \ldots, I_n^i$ are weights associated with a given precision $P$. For each trial a pattern search is generated randomly, but just for representation (excluding order for real process) as $I_x^i = I_1^i + I_2^i + I_3^i \cdots + I_n^i$ and for reaching an interpolated information for $P$, a set of values and trials can be defined as consecutive events for reaching precision $(I_1^i, 1), (I_2^i, 2), (I_3^i, 3), \ldots, (I_n^i, n)$ and a function of $Prob(I^i > I)$ can be defined as,

$$f_i(I^i) = \frac{I^i - I_x^i}{I_n^i} + i \tag{11}$$

And consecutively, the heuristic cognitive linear function searches for environment patterns coding of,

$$f = f_1 I^i, \text{ for any given } P \text{ as } f = f_1 P_1 + f_2 P_2, \ldots, f_n P_n. \tag{12}$$

As organism limitations to process environment patterns fail, it leads to the inverse phenomena of item C.

*D. Discrete probability distribution and monotonically increasing*

**Definition 5.** In the other hand if $I^i < I$, but present definition with limits of $\mathbb{P}((\infty, P]) \leq \mathbb{P}((-\infty, I]) \geq \mathbb{P}((-\infty, I^i])$, and cognitive system presents lower information processing skills, it can be written as,

$$F_p(I) = Prob(P < I) \tag{13}$$

$$\therefore f'(I) - f'(P) = (I - P)f''(I^i) > 0. \text{ Hence, } f'(P) < f'(I).$$

Thus, reducing time and pathways for reaching precision $P$ of information given $I$. Or

$$i = Prob(I^i < I).Prob\left(\frac{T}{I}\right) = Prob\left(\frac{T}{I}\right)^{I^i < I}, \tag{14}$$

*E. Continuous probability density function distribution and monotonically decreasing*

**Definition 6.** Keeping *T* and processing of *I* unstable and unpredictable *a priori*. Following this path, if $I^i$ under high complexity of *n-dimensions*, then *m* assumes distributions like $F(P) > F(I)$. If the probabilistic densities have patterns like $\text{Prob}[P \leq I \leq I^i]$, where,

$$i = Prob(I^i > n).Prob\left(\frac{T}{I}\right) = \int_P^{I^{i>n}} f_I\left(\frac{T}{I}\right) dI^i. \tag{15}$$

*F. Continuous probability density function distribution and monotonically increasing*

**Definition 7.** Otherwise, if *m* express $F(P) < F(I)$ for *n-dimensions*, then,

$$i = Prob(I^i < n).Prob\left(\frac{T}{I}\right) = \frac{dI^i}{dI} F_P(I.I^i < n) \tag{16}$$

*G. Joint probability density function*

**Definition 8.** Considering now *n-dimensions* of external source $(S^i)$, not only individual, *m* assumes distribution like,

$$i = Prob(I^i < S^i).Prob\left(\frac{T}{I}\right), \tag{17}$$

Where $I^i$ assumes notations of,

$$Prob(I^i \to S^i), \tag{18}$$

Giving new equation like,

$$i = Prob S^i.Prob\left(\frac{T}{I}\right) = \int_{S^i}^{I^i} f P_1, \dots, P_n(S^i_1, \dots, S^i_n) dS^i_1 \dots dS^i_n. \tag{19}$$

It is worth mentioning that in time dependent systems, regularity allows the continuous flow of information and possible interruptions caused by the exchange of information between distinct systems generates deceleration of subsequent processes. In other words, the frequency with which activities are performed are dependent on continuous flows to avoid saturation of the work steps that do not have their finalization in the appropriate time. In large information flows, PDFs (probability density function) can be generated on account of chaotic profiles between time-controlled systems [19,20].

## IV. CONCLUSION

This definition, as a rule, is only theoretical and does not prove empirically for any agent or modes of production of a firm when trying to reach always 100% perfect precision, except for artificial intelligence. However, this example illustrates the ideal way of analyzing a workflow for hybrid organizational systems, where the production and precision remains constantly modifying for qualitative and/or quantitative parameters in the sequence of events agent/information/processing/time sequence, which generates precision ideally for a single, a chain or web of events whose expression can be regulated as a mathematical modelling for any region of the system.

It is highly suggested for other researchers in this field or related, to keep in the search of workflows modelling dynamics with new properties of biological systems as an adaptive integration of multidimensional analysis and mathematical approaches.

This article was exclusively performed for qualitative mathematical purposes evaluation and not an in depth calculation analysis for workflows. Qualitative analysis can be applied to flowcharts or other qualitative methods of workflow and information processing analysis. Thou, calculations can be useful for digital workflows and information processing considering for it, several approaches to be implemented by using equation (1) and the other derivations. A very similar approach towards a quantitative method (full mathematical descriptive formulation) for objects to objects interaction and formal framework for probabilistic unclean databases (PUD) was proposed for data analysis in Ilyas *et al* [21]. A full description for *structured predictions* to achieve probabilistic modelling was developed [21] and it is suggested as a tool of information processing as authors require for further investigations of this new techniques, programming parameters that considers biological empirical variances that can serve as input data and more descriptive performance of complex adaptive systems as this issue (quantitative descriptive methods) was not covered by this research.

Corresponding author: charlestelles@seed.pr.gov.br / charles.telles@uol.com.br

Affiliation: Secretary of State for Education and Sport. Curitiba, Paraná, Brazil.
Department: Operational Efficiency and Planning Group.

V. REFERENCES

[1] Oughton EJ, Usher W, Tyler P, Hall JW. Infrastructure as a complex adaptive system. Complexity. 2018;2018.
[2] Cardoso J. Approaches to compute workflow complexity. In Dagstuhl Seminar Proceedings 2006. Schloss Dagstuhl-Leibniz-Zentrum für Informatik.
[3] Dooley KJ. A complex adaptive systems model of organization change. Nonlinear dynamics, psychology, and life sciences. 1997 Jan 1;1(1):69-97.
[4] Miller JH, Page SE. Complex adaptive systems: An introduction to computational models of social life. Princeton university

press; 2009 Nov 28.
[5] Pentland BT, Feldman MS. Organizational routines as a unit of analysis. Industrial and corporate change. 2005 Aug 26;14(5):793-815.
[6] Gabry J, Simpson D, Vehtari A, Betancourt M, Gelman A. Visualization in Bayesian workflow. Journal of the Royal Statistical Society: Series A (Statistics in Society). 2019 Feb;182(2):389-402.
[7] Yang C, Nghiem L, Erdle J, Moinfar A, Fedutenko E, Li H, Mirzabozorg A, Card C. An efficient and practical workflow for probabilistic forecasting of brown fields constrained by historical data. In SPE Annual Technical Conference and Exhibition 2015 Sep 28. Society of Petroleum Engineers.
[8] Lindner B, Auret L, Bauer M. A Systematic Workflow for Oscillation Diagnosis Using Transfer Entropy. IEEE Transactions on Control Systems Technology. 2019 Feb 27.
[9] Stroiteleva TG, Vukovich GG. Mathematical modeling of workflows in production systems. Modern Applied Science. 2015 Mar 1;9(3):173.
[10] Gagnepain P, Ivaldi M. Stochastic frontiers and asymmetric information models. Journal of Productivity Analysis. 2002 Sep 1;18(2):145-59.
[11] Narman P, Buschle M, Konig J, Johnson P. Hybrid probabilistic relational models for system quality analysis. In2010 14th IEEE International Enterprise Distributed Object Computing Conference 2010 Oct 25 (pp. 57-66). IEEE.
[12] Gray RM, Gray RM. Probability, random processes, and ergodic properties. New York: Springer-Verlag; 1988.
[13] Bertsekas DP, Tsitsiklis JN. Introduction to probability. Belmont, MA: Athena Scientific; 2002 Jun.
[14] Sproston J. Decidable model checking of probabilistic hybrid automata. In International Symposium on Formal Techniques in Real-Time and Fault-Tolerant Systems 2000 Sep 20 (pp. 31-45). Springer, Berlin, Heidelberg.
[15] Hofbaur MW, Williams BC. Mode estimation of probabilistic hybrid systems. InInternational Workshop on Hybrid Systems: Computation and Control 2002 Mar 25 (pp. 253-266). Springer, Berlin, Heidelberg.
[16] St-Aubin R, Friedman J, Mackworth AK. A formal mathematical framework for modeling probabilistic hybrid systems. Annals of Mathematics and Artificial Intelligence. 2006 Aug 1;47(3-4):397-425.
[17] Bates CJ, Lerch RA, Sims CR, Jacobs RA. Adaptive allocation of human visual working memory capacity during statistical and categorical learning. Journal of vision. 2019 Feb 1;19(2):11-.
[18] Trübutschek D, Marti S, Dehaene S. Temporal-order information can be maintained in non-conscious working memory. Scientific reports. 2019 Apr 24;9(1):6484.
[19] Telles CR. Geometrical Information Flow Regulated by Time Lengths: An Initial Approach. Symmetry. 2018 Nov 16;10(11):645.
[20] Sheng L, Yushun F, Huiping L. Dwelling time probability density distribution of instances in a workflow model. Computers & Industrial Engineering. 2009 Oct 1;57(3):874-9.
[21] De Sa C, Ilyas IF, Kimelfeld B, Ré C, Rekatsinas T. A Formal Framework for Probabilistic Unclean Databases. In22nd International Conference on Database Theory 2019 Mar.